\documentclass[12pt]{article}
\title{The Poisson Rain Tessellation Under Spatial Expansion and Temporal Transformation}
\author{Eike Biehler, \\Friedrich-Schiller-Universit{\"a}t Jena, \\eikebiehler@web.de}
\date{\today}
\usepackage[latin1]{inputenc}
\usepackage{makeidx}
\usepackage{graphicx}
\usepackage{amssymb}
\usepackage{amsmath}
\usepackage{multirow}
\usepackage[paper=a4paper,left=30mm,right=30mm,top=15mm,bottom=25mm]{geometry}
\usepackage{lscape}
\begin{document}
\newtheorem{Definition}{Definition}
\newtheorem{Satz}{Theorem}
\newtheorem{Lemma}{Lemma}
\newtheorem{Korollar}{Corollary}
\newtheorem{Bemerkung}{Remark}
\newtheorem{Vermutung}{Conjecture}
\newtheorem{Bedingung}{Condition}

\maketitle

\begin{abstract}
For cell-division processes in a window, Cowan introduced four \textit{selection rules} and two \textit{division rules} each of which stands for one cell-division model. One of these is the area-weighted in-cell model. In this model, each cell is selected for division with a probability that corresponds to the ratio between its area and the area of the whole window. This selected cell is then divided by throwing a (uniformly distributed) point into the cell and drawing a line segment through the point under a random angle with the segment ending at the cell's boundary.\\
For the STIT model which uses both a different selection and a different division rule, Martínez and Nagel showed that for a STIT process $\{Y(t, W): t \geq 0\}$ the process $\{e^t Y(e^t, W): t \in \mathbb{R}\}$ is not only spatially stationary but also temporally.\\
For a continuous-time area-weighted in-cell model it is shown by using the different and generalizing approach of a Poisson Rain that in order to get temporal stationarity it is necessary to have an exponential spatial expansion. For the temporal transformation it is found that there is a strong relation between that transformation and the intensity of the Poisson Rain in time.
\end{abstract}

\section{Introduction}
In \cite{Cowan}, Cowan examined cell-division processes in a window that are subdivided according to the ways a cell is chosen for division (by a selection rule) and then actually divided (obeying a division rule). Cowan considered area-, perimeter-, corner-number-weighted and equally-likely selection rules, where e.g. in the area-weighted case the probability for a cell to be chosen for division corresponds to the ratio between its area and the area of the whole window. The two division rules correspond to throwing a point into the cell or onto the cell's boundary (with the point uniformly distributed on the set where it is supposed to fall). It each case, a line is then drawn through that point under a random angle ending at the boundary of the cell.\\
While Cowan examined his processes in discrete time it is not difficult to transform the discrete-time process into a continuous-time process by assigning a lifetime to each cell which is exponentially-distributed according to the weighting parameter of the cell (e.g. the lifetime of a cell $C$ in the area-weighted case is $\mathcal{E}(\lambda_2(C))$-distributed with $\lambda_2$ describing the area).\\
The STIT model which was introduced by Mecke, Nagel and Weiß and whose main characteristics are described in \cite{Nagel:2005} fits into this categorization as an on-boundary Cowan model weighted according to the (generalized) perimeter. For a STIT process $\{Y(t,W): t \geq 0\}$, Martínez and Nagel showed in \cite{Martinez-Nagel} that the process $\{e^t Y(e^t, W): t \in \mathbb{R}\}$ is not only spatially stationary (as any STIT process) but also temporally.\\
It is the goal of this paper to deduce necessary conditions for such a stability property to hold for the continuous-time area-weighted in-cell Cowan model.\\
In order to do this, we first introduce a Poisson Rain with an intensity that may (or may not) vary over the time. It is quite straightforward to see that this Poisson Rain induces such a continuous-time area-weighted in-cell Cowan model as it does not matter whether we first choose a cell area-weightedly and then throw a uniformly-distributed point into the cell or whether we throw a uniformly-distributed point into the whole window and then consider the cell into which the point falls as the chosen cell. Of course, the lifetime distribution of a cell depends on the intensity of the Poisson Rain.\\
To this end, we first examine the properties of the Poisson Rain and what happens when we expand the window in which it comes down and what happens when we transform the time. Then, conditions are deduced for the transformed Poisson Rain to have certain stability properties. Finally, it is explained how to come from the Poisson Rain to a tessellation; this tessellation in a sense inherits the stability conditions of the Poisson Rain and is called SAWSER if the stability conditions are fulfilled.\\
It turns out that a SAWSER tessellation can only exist if the spatial expansion is exponential in time and if the Poisson Rain intensity and the time transformation fulfil a common condition.\\
The idea for this paper came from a discussion with Richard Cowan, Werner Nagel and Viola Weiß on July 1, 2011.


\section{The Poisson Rain}
Let us first, independently of any tessellation, consider a Poisson point process.\\
Let $\chi: \mathbb{R}^2 \times [0, \infty) \rightarrow [0, \infty)$ be a differentiable function. Let further $\tilde{\Phi}_{\chi}$ be a marked Poisson point process on $\mathbb{R}^2 \times [0, \infty) \times [0, \pi)$ with an intensity measure $\chi(x, u) \cdot d (\lambda_2 \otimes \lambda_1 \otimes \mathcal{P})(x, u, \alpha)$ such that for any Borel subsets $S \subset \mathbb{R}^2$, $T \subset [0, \infty)$ and $A \subset [0, \pi)$ the number $\tilde{\Phi}_{\chi}(S \times T \times A)$ of points falling into $S \times T$ with an angle mark from $A$ is determined by
\begin{equation}\label{eq: Definition-Phi-tilde-chi}\begin{array}{l} 
\mathbb{P}(\tilde{\Phi}_{\chi}(S \times T \times A) = k) = \\\\ \frac{1}{k!} \left(\int_{S \times T \times A} \chi(x, u) d(\lambda_2 \otimes \lambda_1 \otimes \mathcal{P})(x, u, \alpha) \right)^k \exp\left(-\int_{S \times T \times A} \chi(x, u) d(\lambda_2 \otimes \lambda_1 \otimes \mathcal{P})(x, u, \alpha) \right)\end{array}\end{equation}
 with $\lambda_1$ on the interval $[0, \infty)$ and $\lambda_2$ on $\mathbb{R}^2$ being the one- and two-dimensional Lebesgue measure and $\mathcal{P}$ being a probability measure, respectively.\\
Each point in $\tilde{\Phi}_\chi$ can be identified as $(x_i, \tau_i, \alpha_i)$ with $\tau_i$ called the time of the point falling onto the plane and $\alpha_i$ being the angle orthogonal to the line to be drawn through $x_i$.\\
We will not need this angle in the next few sections but we will use it when we describe how we come to a tessellation from this point process. Therefore, we will make use of the (unmarked) Poisson point process $\Phi_\chi$ which we can deduce from $\tilde{\Phi}_\chi$ via the projection \begin{equation}\label{eq: Projection-Phi-chi}\Phi_\chi(S \times T) = \tilde{\Phi}_\chi(S \times T \times [0, \pi)).\end{equation}
The mean number of points of $\Phi_\chi$ falling into $S$ during $T$ is
$$\mathbb{E}\Phi_\chi(S \times T) = \int_{S \times T} \chi(x, u) d(\lambda_2 \otimes \lambda_1)(x, u).$$
This can also be written as
\begin{equation}\label{eq: Poisson-Rain-general}\mathbb{E}\Phi_\chi(S \times T) = \int_{\mathbb{R}^2 \times [0, \infty)} \mathbf{1}_{S \times T}(x, u) \chi(x, u) d(\lambda_2 \otimes \lambda_1)(x, u)\end{equation} with $\mathbf{1}_D$ being the indicator function for a set $D$.\\
With $\chi \equiv 1$, we get the usual homogeneous Poisson point process in space \textit{and} time with intensity $1$.\\
We can interpret this Poisson point process as a \textit{Poisson Rain} on the plane with the rain drops (the points $x_i$) touching the plane at the  time $\tau_i$. The intensity of this rain in space and time is governed by the function $\chi$. In the following, we will consider the intensity of the rain as independent of the location, therefore we will write only $\chi(u)$. Thus, we will now use the intensity measure $\chi(u) \cdot d (\lambda_2 \otimes \lambda_1)(x, u)$.

\section{Two approaches for expansion into the plane}
\label{sec: Two-approaches}
We now consider a convex and compact set $S \subset \mathbb{R}^2$. We want this set to grow with the time and want to consider the mean number of points falling into the \textit{growing} subset. Here, we have two possibilities to describe this growth until an arbitrary time $t$:
\begin{itemize}
	\item The first approach, to be called Growth During Rain (GDR), will see the set grow at the same time as the points are falling, so the growth of the set and the falling of points in the Poisson Rain are considered to take place simultaneously.
	\item The second approach, to be called Growth After Rain (GAR), will see first the points falling into a non-growing set until $t$ and \textit{then} grow the set to its size at time $t$.
\end{itemize}
In both approaches, a point that has fallen into the set moves away from the origin as the time (and the set) grows.\\
It should be clear that the GDR approach yields a stronger connection between the growth of the window and number of points falling into it than the GAR approach.\\
Let us use monotonically increasing, differentiable functions $\xi_D: [0, \infty) \rightarrow [0, \infty)$ and $\xi_A: [0, \infty) \rightarrow [0, \infty)$ for the GDR and GAR approach, respectively, to describe the growth of the set. Let us have a reference time $s_0^D$ and $s_0^A$ for which $\xi_D(s_0^D) = 1$ and $\xi_A(s_0^A) = 1$ holds. Of course, this reference time may depend on the function $\xi_D$ or $\xi_A$, respectively. Accordingly, by $\chi_D$ and $\chi_A$ we describe the intensities of the rain in the corresponding approach.\\
Usually the time set $T$ will now be an interval $[s, s+t]$ with $s \geq 0$ and $t \geq 0$.

\subsection{The GDR approach}
For the function $\xi_D$, let us define $S = S_{s_0^D} = \xi_D(s_0^D) S$ which is supposed to be a compact and convex set in the plane (a window). Thus we consider $S$ to be our reference set with its size at the time $s_0^D$. Then we can define $S_s = \xi_D(s) S$ and $S_{s+t} = \xi_D(s+t) S$. For any time $u \in [s, s+t]$, we can define $S_u = \xi_D(u) S$. Such a set has the area $$\lambda_2(S_u) = \lambda_2(\xi_D(u) S) = \xi_D^2(u) \lambda_2(S).$$ Also, $$S_u = \xi_D(u) S = \frac{\xi_D(u)}{\xi_D(s)} S_s.$$
Thus, we can express the set $S_u$ (whose area grows with growing $u$) as a function of the set $S_s$ at the time $s$. A point $(x_i, \tau_i)$ of the Poisson Rain $\Phi_{\chi_D}$ that hits the plane at the time $\tau_i$ hits the set $S_{\tau_i}$ if and only if $x_i \in S_{\tau_i}$. At the time $s+t$, this point has moved from $x_i$ to $\frac{\xi_D(s+t)}{\xi_D(\tau_i)} x_i$. Thus, a point that fell at the time $\tau_i$ is in $S_{s+t}$ if and only if it fell into $S_{\tau_i}$.
With the point process $\Phi_{\chi_D}$ as defined for $\chi$ through equations (\ref{eq: Projection-Phi-chi}) and (\ref{eq: Definition-Phi-tilde-chi}) in hand, we can define for any time $s$ with a set $S_s$, a Borel set $A \subset \mathbb{R}^2$ and an interval $[s, s+t]$ the point process $$\Phi_{\chi_D, \xi_D, s, S_s; s+t} = \left\{\left(\frac{\xi_D(s+t)}{\xi_D(\tau_i)}x_i, \tau_i\right) : (x_i, \tau_i) \in \Phi_{\chi_D}, x_i \in S_{\tau_i}, \tau_i \in [s, s+t]\right\}$$ which depends on the functions $\chi_D$ and $\xi_D$,  on the time $s$ and the examined set $S_s$ as well as on the time $s+t$ we consider.\\Thus, the mean number of points falling into that growing window $S_u = \xi_D(u) S$ for $u \in [s, s+t]$ can be calculated in analogy to (\ref{eq: Poisson-Rain-general}) and considering the point process as a measure now as
$$\begin{array}{rl}& \mathbb{E}\Phi_{\chi_d, \xi_D, s, S_s; s+t}(S_{s+t} \times [s, s+t])\\&\\
 = & \int_{\mathbb{R}^2 \times [0, \infty)} \mathbf{1}_{S_u \times [s, s+t]}(x, u) {\chi_D}(u) d(\lambda_2 \otimes \lambda_1)(x, u)\\&\\
 = & \int_{[0, \infty)} \mathbf{1}_{[s, s+t]}(u) \left(\int_{\mathbb{R}^2} \mathbf{1}_{S_u}(x) d\lambda_2(x)\right) {\chi_D}(u) d\lambda_1(u)\\&\\
 = & \int_s^{s+t} \lambda_2(S_u) {\chi_D}(u) du\\&\\
 = & \int_s^{s+t} \frac{\xi_D^2(u)}{\xi_D^2(s)} \lambda_2(S_s) {\chi_D}(u)du\\&\\
 = & \frac{\lambda_2(S_s)}{\xi_D^2(s)} \int_s^{s+t} \xi_D^2(u) {\chi_D}(u)du. \end{array}$$
All points of $\Phi_{\chi_D, \xi_D, s, S_s; s+t}$ are per definition contained within $S_{s+t}$.\\ 
This process is a Poisson point process: The points have moved away from the origin at the same speed as the window grows. This conserves the property that the mean number of points within a set is proportional to the area of the set. Additionally, all points have fallen independently of each other.

\subsection{The GAR approach}
Let us define again that $S = S_{s_0^A} = \xi_A(s_0^A) S$ as a compact and convex set in the plane. Then we can say again that $S_s = \xi_A(s) S$ and $S_{s+t} = \xi_A(s+t) S$.\\
From the point process that is now denoted $\Phi_{\chi_A}$ we can define a point process $$\Phi_{\chi_A, \xi_A, s, S; s+t} = \left\{\left(\xi_A(s+t) x_i, \tau_i\right): (x_i, \tau_i) \in \Phi_{\chi_A}, x_i \in S, \tau_i \in [s, s+t]\right\}.$$ Note that the reference here is made to the set $S$ which exists at the reference time $s_0^A$ with $\xi_A(s_0^A)=1$.\\
We then get 
$$\begin{array}{rl} & \mathbb{E}\Phi_{{\chi_A}, \xi_A, s, S; s+t}(S_{s+t} \times [s, s+t])\\&\\
= & \mathbb{E}\Phi_{\chi_A}(S \times [s, s+t])\\&\\
= & \mathbb{E}\Phi_{\chi_A}\left(\frac{1}{\xi_A(s)}S_s \times [s, s+t]\right)\\&\\
= & \lambda_2\left(\frac{1}{\xi_A(s)} S_s \right) \int_s^{s+t} {\chi_A}(u) du\\&\\
= & \frac{\lambda_2(S_s)}{\xi_A^2(s)} \int_s^{s+t} {\chi_A}(u)du.\end{array}$$
The only difference to the GDR method is that the function $\xi_A^2(u)$ does \textit{not} appear in the integral in contrast to $\xi_D^2(u)$ in the GDR case.\\
This process is a Poisson point process as well: The point process $\Phi_{\chi_A}$ is a Poisson process which is only stretched to get the point process in question, thus keeping the Poisson properties of the unstretched process.

\section{Time transformation}
\label{sec: Time-transformation}
\subsection{The function $\psi$ for both approaches}
We can now also transform the time. Let us have a function $\psi: (-\infty, \infty) \rightarrow [0, \infty)$ which is differentiable on the support $\{s \in \mathbb{R}: \psi(s) > 0\}$ of $\psi$.\\
For the GDR method, by the $\xi_D^2$ in the integral we have an intricate link between the number of points falling into the \textit{growing} area and the area of the set examined. Therefore, the only process that is a time-transformed GDR process is
$\Phi_{{\chi_D}, \xi_D,  \psi; s, S_s, s+t}$ which is defined by
$$\Phi_{{\chi_D}, \xi_D,  \psi; s, S_s, s+t}(S_{s+t} \times [s, s+t]) = \Phi_{\chi_D, \xi_D, \psi(s), S_{\psi(s)}; \psi(s+t)}(S_{\psi(s+t)} \times [\psi(s), \psi(s+t)])$$ for $s \in \mathbb{R}$ and $t \geq 0$. We then get
\begin{equation}\label{eq: GDR-method-points-falling-general}\begin{array}{rl}& \mathbb{E}\Phi_{{\chi_D}, \xi_D,  \psi; s, S_s, s+t}(S_{s+t} \times [s, s+t])\\&\\ = & \mathbb{E}\Phi_{\chi_D, \xi_D, \psi(s), S_{\psi(s)}; \psi(s+t)}(S_{\psi(s+t)} \times [\psi(s), \psi(s+t)])\\&\\ = & \frac{\lambda_2(S_{\psi(s)})}{\xi_D^2(\psi(s))} \int_{\psi(s)}^{\psi(s+t)} \xi_D^2(u) {\chi_D}(u) du.\end{array}\end{equation}
For the GAR method, we have a fixed window $S$ into which during the time interval $[\psi(s), \psi(s+t)]$ a certain number of points fall. \textit{Afterwards}, we let the window grow; in this growth, the \textit{number} of points in the growing set remains constant but the points themselves \textit{move} away from the origin. Thus, we can set $\xi$ independent of $\psi$ and define $\Phi_{{\chi_A}, \xi_A, \psi; s, S, s+t}$ as $$\begin{array}{rl} & \Phi_{{\chi_A}, \xi_A, \psi; s, S, s+t}\\&\\ = & \left\{(\xi_A(s+t) x_i, \tau_i): (x_i, \tau_i) \in \Phi_{\chi_A}, x_i \in S, \tau_i \in [\psi(s), \psi(s+t)]\right\}\\&\\ = &  \Phi_{{\chi_A}, \xi_A, \psi(s), S; \psi(s+t)}\end{array}.$$
We then get \begin{equation}\label{eq: General-definition-GAR}\begin{array}{rl}&\mathbb{E}\Phi_{{\chi_A}, \xi_A, \psi; s, S, s+t}(S_{s+t}\times [s, s+t]) \\&\\= & \mathbb{E} \Phi_{{\chi_A}, \xi_A, \psi(s), S; \psi(s+t)} (S_{s+t} \times [\psi(s), \psi(s+t)]) \\&\\
=  & \lambda_2(S) \int_{\psi(s)}^{\psi(s+t)} {\chi_A}(u) du\\&\\
=  & \frac{\lambda_2(S_s)}{\xi_A^2(s)} \int_{\psi(s)}^{\psi(s+t)} {\chi_A}(u) du.\end{array}\end{equation}

\subsection{Equivalence}
If we start our consideration at an arbitrary time $s$ (which is then transformed by $\psi$) and make the condition that the sets we are considering at that time ($S_{\psi(s)}$ in the GDR, $\xi_A(s) S$ in the GAR approach as $\psi$ does only have an influence on the number of points within the time interval but not the size of the set) have the same area, we can see from equations (\ref{eq: GDR-method-points-falling-general}) and (\ref{eq: General-definition-GAR}) that with \begin{equation}\label{eq: Equivalence-1}(\xi^2_D \circ \psi)(s) = \xi^2_A(s)\end{equation} for any $s \in \mathbb{R}$ and \begin{equation}\label{eq: Equivalence-2}(\xi^2_D \cdot \chi_D)(u) = \chi_A(u)\end{equation} for any $u \in [0, \infty)$ we have the same mean number of points falling into that set $S_{\psi(s+t)}$ for the GDR approach and $S_{s+t}$ for the GAR approach, respectively, until the time $s+t$.\\
Note that from equation (\ref{eq: Equivalence-1}) the requirement $\psi(s_0^A) = s_0^D$ follows immediately.

\subsection{Conditions for the functions $\chi$, $\xi$ and $\psi$}
Due to the equivalence of the approaches as per equations (\ref{eq: Equivalence-1}) and (\ref{eq: Equivalence-2}), in the following, we will predominantly consider the GAR approach. In doing so, we will not always use the indices $_A$ and $_D$ when it is not necessary to make a distinction between the two approaches.\\
In the first sections, we have stated that the functions $\chi$, $\xi$ and $\psi$ shall be differentiable and non-negative. The requirement of differentiability will later make some calculations easier. Additionally, we want $\psi$ to be monotonically increasing as it describes the time span we consider, thus $\psi(s+t) \geq \psi(s)$ should be true for $t \geq 0$. It should furthermore not be constant within the support $\{s \in \mathbb{R}: \psi(s) > 0\}$ as then, every interval lying within such an area of constancy becomes a point yielding any integral zero. Thus, $\psi$ is supposed to be strictly monotonically increasing on its support.
Additionally, $\xi^2$ should be monotonically increasing as otherwise speaking of a growing window would not make any sense; there also should be a time $s_0 \in [0, \infty)$ such that $\xi(s_0)=1$. When transforming the time when we want also to have an $s \in (-\infty, \infty)$ such that $\psi(s) = s_0$ and thus $\xi(\psi(s)) = 1$.\\
This leads to the following conditions:
\begin{Bedingung}
\label{Bedingung: Regularity}
We only admit non-negative functions $\chi$, $\psi$ and $\xi$ that are differentiable and for which $$\frac{d}{dt} \psi(t) > 0 \textrm{ for any $t \in \mathbb{R}$ with $\psi(t) > 0$}$$ and $$\frac{d}{dt} \xi^2(t) \geq 0$$ holds. There also must be a time $s_0 \in [0, \infty)$ such that $\xi(s_0)=1$ and a time $s \in (-\infty, \infty)$ such that $\psi(s)=s_0$.
\end{Bedingung}

\section{Stability considerations}
\label{sec: Stability-considerations}
For STIT tessellation processes, Martínez and Nagel showed in \cite{Martinez-Nagel} that for a STIT tessellation process $(Y_t, t \geq 0)$ the transformed process $(e^t Y_{e^t}, t \in \mathbb{R})$ is stationary. While we have not yet shown how a tessellation can arise from any of our two approaches, we can motivate two conditions for stability of our point process.
\begin{itemize}
	\item The \textbf{first condition} for stability is that the mean number of points per unit area becomes constant (and greater than zero) at least asymptotically. Otherwise, the process may degenerate by either having ever less or ever more points contained within a set or it may oscillate around a certain value without actually approaching a limit. We will call that the MEPA condition below.
	\item The \textbf{second condition} is that, given a set at an arbitrary time $s$, the waiting time for this (growing and moving) set to be hit does not depend on the time $s$. We will call this the INOT condition below.
\end{itemize}

\subsection{The MEPA condition}
The \underline{me}an number of \underline{p}oints per unit \underline{a}rea (MEPA) is the ratio of the mean number of points in a set $S_{s+t}$ at a time $s+t$ and the area of that set $S_{s+t}$. Thus, we can write for some $s \in \mathbb{R}$
$$\begin{array}{rl}M = & \lim_{t \rightarrow \infty} \frac{\mathbb{E}\Phi_{{\chi_A}, \xi_A, \psi; s, S, s+t}(S_{s+t} \times [s, s+t])}{\lambda_2(S_{s+t})}\\&\\
= & \lim_{t \rightarrow \infty} \frac{\xi_A^2(s)}{\xi_A^2(s+t) \lambda_2(S_s)} \cdot \frac{\lambda_2(S_s)}{\xi_A^2(s)} \int_{\psi(s)}^{\psi(s+t)} {\chi_A}(u) du \\&\\
= & \lim_{t \rightarrow \infty} \frac{\int_{\psi(s)}^{\psi(s+t)} {\chi_A}(u) du}{\xi_A^2(s+t)}
\end{array}$$ if that limit exists.
This leads to 
\begin{Bedingung}
\label{Bedingung: MEPA}
Let the limit $\lim_{t \rightarrow \infty} \frac{1}{\xi_A^2(s+t)} \int_{\psi(s)}^{\psi(s+t)} {\chi_A}(u) du$ exist. Then, the inequality
$$0 < \lim_{t \rightarrow \infty} \frac{1}{\xi_A^2(s+t)} \int_{\psi(s)}^{\psi(s+t)} {\chi_A}(u) du = M < \infty$$ must hold for every $s \in \mathbb{R}$.
\end{Bedingung}

\subsection{The INOT condition}
The \underline{in}dependence of the \underline{o}bservation \underline{t}ime (INOT) condition is to ensure that whenever we observe a (growing and moving) set the waiting time for a point of the point process to fall into it has the same distribution.\\
Let us consider a set $S_s \subset \mathbb{R}^2$ at the time $s$.  Then, the distribution of the waiting time $X_{S_s}$ for $S_s$ to be hit after the commencement (at the time $s$) of the observation can be described as
\begin{equation}\label{eq: waiting-time-INOT}\mathbb{P}(X_{S_s} > t) = \mathbb{P}(\Phi_{{\chi_A}, \xi_A, \psi; s, S, s+t}(S_{s+t}\times [s, s+t]) = 0) = \exp\left(- \frac{\lambda_2(S_s)}{\xi_A^2(s)} \int_{\psi(s)}^{\psi(s+t)}  \chi_A(u) du\right),\end{equation} respectively.
It is this waiting time's distribution which must be independent of the time $s$. Let us have a time $s_1$ and a set $S^{(1)}$ on one hand and another time $s_2$ and another set $S^{(2)}$ on the other hand. Let the two sets have, however, equal area, i.e. $\lambda_2(S^{(1)}) = \lambda_2(S^{(2)})$. As $\lambda_2(S^{(1)}) = \lambda_2(\xi_A(s_1) S_1) = \xi_A^2(s_1) \lambda_2(S_1)$ for some set $S_1$ at the reference time $s_0^A$ and, analogously, $\lambda_2(S^{(2)}) = \lambda_2(\xi_A(s_2) S_2) = \xi_A^2(s_2) \lambda_2(S_2)$ for some set $S_2$ at the reference time $s_0^A$, we can state for the sets $S_1$ and $S_2$ at the reference time $$\lambda_2(S_1) \neq \lambda_2(S_2) \Leftrightarrow \xi_A(s_1) \neq \xi_A(s_2).$$ Because the sets $S^{(1)}$ and $S^{(2)}$ do have equal area however, we can state that the term $$\frac{1}{\xi_A^2(s)} \int_{\psi(s)}^{\psi(s+t)} \chi_A(u) du$$ must be independent of $s$.
Therefrom, we can deduce
\begin{Bedingung}
\label{Bedingung: INOT}
The equation
$$\frac{d}{ds} \left(\frac{1}{\xi_A^2(s)} \int_{\psi(s)}^{\psi(s+t)} \chi_A(u) du\right) = 0$$ must hold for every $t > 0$ and every $s \in \mathbb{R}$.
\end{Bedingung}
Note that as all functions are assumed differentiable the whole term is differentiable as well.

\section{The SAWSER tessellation}
\label{sec: SAWSER}
We will first explain how we get a tessellation $Y$ from the Poisson point process $\tilde{\Phi}_{\chi}$. Then, the terms of a tessellation with an area-weighted selection rule (AWSER) and of a stable AWSER tessellation (SAWSER) will be defined. Lastly, we will prove that the tessellation $Z$ arising from spatial expansion and temporal transformation of the original tessellation $Y$ is indeed a SAWSER tessellation.

\subsection{The unexpanded and untransformed tessellation in $W$}
As the number of points falling into a bounded spatial set $W$ during a bounded time interval $[0, t]$ is almost surely finite, one can order these points according to the time they fall. Let $\kappa = \tilde{\Phi}_{\chi}(W \times [0, t] \times [0, \pi))$ be the (random) number of points falling into $W$ during $[0, t]$ with any angle. Then we can with probability 1 re-arrange the indices in $\tilde{\Phi}_{\chi}(W \times [0,t] \times [0, \pi)) = \{(x_1, \tau_1, \alpha_1), ..., (x_\kappa, \tau_\kappa, \alpha_\kappa)\}$ such that $0 < \tau_1 < ... < \tau_\kappa < t$.\\
At the time $\tau_1$, the point $x_1$ falls into $W$ and a line segment $\gamma_1(x_1, \alpha_1) \cap W$ is drawn through that point with the angle $\alpha_1$. Let us denote $$Y(\tau_1, W) = \gamma_1(x_1, \alpha_1) \cap W$$ as the tessellation in $W$.
Let us now define $$C(x_2, Y(\tau_1, W)) = cl(\{y \in W: \overline{x_2 y} \cap Y(\tau_1, W) = \emptyset\})$$ as the cell of $Y(\tau_1, W)$ in which $x_2$ falls. Here, $cl(S)$ means the closure of a set S and $\overline{x_2 y}$ is the connecting linear segment between $x_2$ and $y$.\\
Then, we define $$Y(\tau_2, W) = Y(\tau_1, W) \cup (\gamma_1(x_2, \alpha_2) \cap C(x_2, Y(\tau_1, W)))$$ as the tessellation after the second point has fallen into the window.\\
Thus, we iteratively define $$C(x_n, Y(\tau_{n-1}, W) = cl(\{y \in W: \overline{x_n y} \cap Y(\tau_{n-1}, W) = \emptyset\})$$ and
$$Y(\tau_n, W) = Y(\tau_{n-1}, W) \cup (\gamma_1(x_n, \alpha_n) \cap C(x_n, Y(\tau_{n-1}, W))).$$
Finally, under the given situation with $\kappa$ points in $W$ during $[0, t]$ we have $$Y(t, W) = Y(\tau_\kappa, W).$$

\subsection{AWSER and SAWSER tessellations}
Let us now consider functions $\xi$ and $\psi$ as defined above. Then we define
$$Z(t, W_t) = \xi(t) Y(\psi(t), W)$$ as the tessellation arising from $Y(\psi(t),W)$ in a window $W_t = \xi(t) W$ and $(Z(t, W_t): t \geq 0)$ as the tessellation process.\\
Let us now define in general the following.
\begin{Definition}
\label{Definition: AWSER}
A tessellation $Z = Z(t, W_t)$ in a compact and convex window $W_t = \xi(t) W \subset \mathbb{R}^2$ with non-empty interior is called AWSER if it employs an area-weighted rule for the selection of cells to be divided at each division step, i.e. if there is a cell division at a certain time $t$ in the window $W_t$, the probability for a certain cell $C \in Cells(Z(t, W_t))$ to be divided by the dividing segment $S_t$ is $$\mathbb{P}(C \cap S_t \neq \emptyset|W_t \cap S_t \neq \emptyset) = \frac{\lambda_2(C)}{\lambda_2(W_t)}.$$
\end{Definition}

\begin{Definition}
\label{Definition: SAWSER}
An AWSER tessellation is called \underline{s}table (and then a \underline{S}AWSER tessellation) if 
\begin{itemize}
\item the ratio of the mean number of cells $\mathbb{E}|Cells(Z(t,W))|$ and the area $\lambda_2(W_t)$ of the window converges towards a positive constant for $t \rightarrow \infty$: $$0 < \lim_{t \rightarrow \infty} \frac{\mathbb{E}|Cells(Z(t,W_t))|}{\lambda_2(W_t)} = M_S < \infty;$$
	\item for each cell $C$, its lifetime $X_C$ is independent of the time of observation, i.e. if there are two cells $C_1$ and $C_2$ of the same area ($\lambda_2(C_1) = \lambda_2(C_2) = L$) at different observation times $s_1$ and $s_2$, then $$\begin{array}{rl} & \mathbb{P}(X_{C_1} > t|C_1 \in Cells(Z(s_1, W_{s_1}), \lambda_2(C_1) = L)\\&\\ = & \mathbb{P}(X_{C_2} > t|C_2 \in Cells(Z(s_2, W_{s_2}), \lambda_2(C_2) = L)\end{array}$$ for any $t \geq 0$ and any $L > 0$.
\end{itemize}
\end{Definition}
Because of the Conditions \ref{Bedingung: MEPA} and \ref{Bedingung: INOT}, we can state
\begin{Lemma}
\label{Lemma: Stability-of-Z}
If the functions $\chi$, $\xi$ and $\psi$ fulfil the Conditions \ref{Bedingung: Regularity} on the regularity, \ref{Bedingung: MEPA} (MEPA) and \ref{Bedingung: INOT} (INOT), then the tessellation $Z(t, W_t)$ is SAWSER.
\end{Lemma}
\textbf{Proof}\\
We have to prove that $Z(t, W_t)$ is AWSER and that it is stable. Let us turn first to the AWSER property.\\
For each point $s$ in time, a cell $C_s$ has a lifetime $X_{C_s}$ for which the probability to exceed a time $t$ can be described per equation (\ref{eq: waiting-time-INOT}) by an exponential function with the cell's area at the observation time in the argument as $$\mathbb{P}(X_{C_s} > t) = \exp\left(- \frac{\lambda_2(C_s)}{\xi_A^2(s)} \int_{\psi(s)}^{\psi(s+t)}  \chi_A(u) du\right).$$ Accordingly, the probability that the waiting time $$X^{W_s} = \min\{X_{C_s}: C_s \in Cells(Z(s, W_s))\}$$ of the state of the whole tessellation in $W_s$ to change exceeds $t$ is $$\mathbb{P}(X^{W_s} > t) = \exp\left(- \frac{\lambda_2(W_s)}{\xi_A^2(s)} \int_{\psi(s)}^{\psi(s+t)}  \chi_A(u) du\right).$$ Due to well-known properties of the exponential distributions, the probability that this cell $C_s$ is the next to be divided is the ratio of its area and the area of the whole window $W_s$. As the area of the cell and of the whole window grow proportionally, this ratio remains the same even within the growing window. This, however, means that the selection rule is area-weighted.\\
The stability follows from the fact that for the number of cells $|Cells(Z(t, W_t))|$ the equation $$|Cells(Z(t, W_t))| = \Phi_\chi(W \times [0, \psi(t)) + 1$$ holds. From Condition \ref{Bedingung: MEPA} we get $$\begin{array}{rl}& \lim_{t \rightarrow \infty}\frac{\mathbb{E}|Cells(Z(t, W_t))|}{\lambda_2(W_t)}\\&\\
= & \lim_{t \rightarrow \infty}\frac{\mathbb{E}\Phi_\chi(W \times [0, \psi(t)) + 1}{\lambda_2(W_t)}\\&\\
= & \lim_{t \rightarrow \infty}\frac{\mathbb{E}\Phi_\chi(W \times [0, \psi(t))}{\lambda_2(W_t)} + \lim_{t \rightarrow \infty}\frac{1}{\lambda_2(W_t)}\\&\\
\stackrel{C2}{=} & M + \lim_{t \rightarrow \infty}\frac{1}{\lambda_2(W_t)}\\&\\
< & \infty\end{array}$$ due to $\frac{d}{dt} \xi^2(s) \geq 0$ for any $s \geq 0$.\\Due to $M > 0$, of course also $\lim_{t \rightarrow \infty}\frac{\mathbb{E}|Cells(Z(t, W_t))|}{\lambda_2(W_t)} > 0$.\\
A cell is a set in the window, so Condition \ref{Bedingung: INOT} makes sure that the waiting time for a point to fall into that cell (i.e. the cell's lifetime) has a distribution which is independent of the time of observation.\\
Thus, $Z(t, W_t)$ is stable and has an area-weighted selection rule; it is a SAWSER tessellation. \hfill $\Box$\\

\section{Dependencies between $\chi$, $\xi$ and $\psi$}
\label{sec: Dependencies-chi-xi-psi}

\subsection{A condition for $\xi^2$}
We have the MEPA and the INOT condition in which we can re-write the integral as
$$\int_{\psi(s)}^{\psi(s+t)} \chi_A(u)du = X(\psi(s+t)) - X(\psi(s)) = (X \circ \psi) (s+t) - (X \circ \psi) (s)$$ with $X$ being the anti-derivative of $\chi_A$.\\
Due to the INOT condition, we can say
$$(X \circ \psi) (s+t) - (X \circ \psi) (s) = f(t) \xi^2_A(s)$$ for some function $f$. As we must separate $s$ and $t$, we need $$(X \circ \psi)(s+t) = (X \circ \psi)(s) \cdot \tilde{f}(t)$$ or $\frac{g(s+t)}{g(s)} = \tilde{f}(t)$ for $g = X \circ \psi$ to hold.\\ Differentiation to $s$ leads to $$0 = \frac{d}{ds}\tilde{f}(t) = \frac{d}{ds} \frac{g(s+t)}{g(s)} = \frac{g'(s+t)g(s) - g(s+t)g'(s)}{g^2(s)},$$ and thus $$g'(s+t)g(s) - g(s+t)g'(s)=0$$ must hold for any $t \geq 0$. Then, however, $$\frac{g'(s)}{g(s)} = \frac{g'(s+t)}{g(s+t)}$$ must hold for any $t \geq 0$. This means that $$\frac{g'(s)}{g(s)} = const.$$ must hold. From the differential equation $$g'(s) = b g(s)$$ we get $$g(s) = a e^{bs} = (X \circ \psi)(s)$$ for $a, b \in \mathbb{R}$.\\
So, we get the equation $$ae^{b(s+t)} - ae^{bs} = f(t) \xi^2_A(s)$$ and thus because of the condition $\xi_A(s_0) = 1$ $$f(t) = a e^{bs_0} \left(e^{bt} - 1\right) \textrm{ and } \xi^2_A(s) = e^{-bs_0}e^{bs}.$$
Putting this into the MEPA condition, we get that
$$0 < \lim_{t \rightarrow \infty} e^{-b(s+t)} a \left(e^{b(s+t)} - e^{bs}\right) = \lim_{t \rightarrow \infty} a\left(1 - e^{-bt}\right) = a = M < \infty.$$ Thus, the MEPA condition is fulfilled as well for $b > 0$.\\
It is noteworthy that $b=0$, i.e. $\xi^2_A(t) = 1$ for all $t \geq 0$, is not admissible as the MEPA condition would not hold then (the limit would be zero as $(X \circ \psi) (s+t) - (X \circ \psi) (s) = 0$).

\subsection{Conditions for $\chi$ and $\psi$}
Per the condition $$(X \circ \psi)(s) = a e^{bs}$$ we get $$\chi(\psi(s)) \frac{d}{ds} \psi(s)= abe^{bs}.$$
As $\frac{d}{ds} \psi(s) > 0$, we can define $\psi^{-1}$ as the inverse function of $\psi$ such that for any $s$ with $\psi(s) > 0$ we have $\psi^{-1}(\psi(s))=s$. We will write $\frac{d}{ds} \psi(s) = \psi'(s)$. Then we have $$\chi(\psi(s)) = \frac{abe^{bs}}{\psi'(s)}$$ and then, inserting $s=\psi^{-1}(t)$,
\begin{equation}\label{eq: chi}\chi(t) = \chi(\psi(\psi^{-1}(t))) = \frac{ab e^{b\psi^{-1}(t)}}{\psi'(\psi^{-1}(t))}.\end{equation}
This is e.g. fulfilled for (Case A) \begin{equation}\label{eq: Case-A}\psi(s) = a e^{bs} \textrm{ and } \chi(s) = 1\end{equation} (due to $\psi^{-1}(s) = \frac{1}{b} \ln \frac{s}{a}$) or for (Case B) $$ \psi(s) = s \mathbf{1}_{[0, \infty)}(s) \textrm { and } \chi(s) = abe^{bs}$$ due to $\psi^{-1}(s) = s = \psi(s)$ on the support of $\psi$. Here, $\mathbf{1}$ is the indicator function which is necessary to have $\psi$ defined on whole $\mathbb{R}$ but map into $[0, \infty)$ only; such a $\psi$ is differentiable on its support as required.\\
For the classification Cowan introduced in \cite{Cowan}, the special case $\chi(s) = 1$ corresponds to the model Cowan-2a, namely the area-weighted model with \textit{in-cell} division rule. As we can deduce from equation (\ref{eq: Case-A}), we \textit{must} transform the time by $\psi(s) = ae^{bs}$ in order to get a stable tessellation.

\subsection{Theorem}
\begin{Satz}
A tessellation $$Z(t, W) = \xi(t) Y(\psi(t), W)$$ with $Y(t, W)$ arising from a Poisson point process $\tilde{\Phi}_\chi$ is SAWSER if and only if the conditions $$\xi^2(t) = e^{-bs_0} e^{bt}$$ for $b > 0$ and $\xi(s_0)=1$ and $$\chi(\psi(s)) \frac{d}{ds} \psi(s)= abe^{bs}$$ are fulfilled.\\
For any given $\psi$ fulfilling the Conditions, $$\chi(t) = \frac{abe^{b \psi^{-1}(t)}}{\psi'(\psi^{-1}(t))}$$ is the fitting rain intensity.
\end{Satz}

\subsection{The re-translation to the GDR process}
For the re-translation to the GDR process, we can use equations (\ref{eq: Equivalence-1}) and (\ref{eq: Equivalence-2}). We have a $\psi$ which only must be strictly monotonically increasing, $\xi_A(t) = e^{-bs_0} e^{bt}$ and $\chi$ as described per equation (\ref{eq: chi}).\\
We can then quite straightforwardly see that $$\xi^2_D(t) = e^{-bs_0} e^{b\psi^{-1}(t)}$$ and $$\chi_D(t) = \frac{ab}{\psi'(\psi^{-1}(t))}.$$
Let us consider two combinations of functions $(\chi, \xi^2, \psi)$, namely $(1, e^{bt}, e^{bt})$ and $(abe^{bt}, e^{bt}, t \mathbf{1}_{[0, \infty)}(t))$, i.e. functions with $s_0=0$. We then get
$$\begin{array}{lrcl}(\textrm{Case A}) & (\chi_A, \xi^2_A, \psi)(t) = (1, e^{bt}, e^{bt}) & \Longleftrightarrow & (\chi_D, \xi^2_D, \psi)(t) = (a t^{-1}, t, e^{bt})\\&&&\\
\textrm{and} &&&\\&&&\\
(\textrm{Case B}) & (\chi_A, \xi^2_A, \psi)(t) = (abe^{bt}, e^{bt}, t \mathbf{1}_{[0, \infty)}(t)) & \Longleftrightarrow & (\chi_D, \xi^2_D, \psi)(t) = (ab, e^{bt}, t \mathbf{1}_{[0, \infty)}(t)).
\end{array}$$
It should be noted that in Case A in the GDR approach $\chi$ is not defined for $t=0$.

\begin{center}
\textbf{Acknowledgment}
\end{center}
I am once more thankful to Werner Nagel for his advice during the development of this article.

\bibliographystyle{plain} \bibliography{literatur}

\end{document}